\documentclass[11pt,a4paper]{article}
\usepackage{tipa}
\usepackage{amsfonts}
\usepackage{amsmath}
\usepackage{amssymb}
\usepackage{amsthm}
\usepackage{latexsym}
\usepackage{amsbsy}
\usepackage[all]{xypic}
\usepackage{amscd}

\newtheorem{thm}{Theorem}[section]
\newtheorem{cor}[thm]{Corollary}
\newtheorem{lem}[thm]{Lemma}

\newtheorem{prop}[thm]{Proposition}
\newtheorem{defn}[thm]{Definition}

\begin{document}

\begin{center}
{\Large \bf Cluster combinatorics of $d-$cluster
categories\footnote{Supported by the NSF of China (Grant 10771112)
and in part
 by Doctoral
Program Foundation of Institute of Higher Education(2009) }}
\bigskip

{\large Yu Zhou and
 Bin Zhu}
\bigskip

{\small
\begin{tabular}{cc}
Department of Mathematical Sciences & Department of Mathematical
Sciences
\\
Tsinghua University & Tsinghua University
\\
  100084 Beijing, P. R. China &   100084 Beijing, P. R. China
\\
{\footnotesize E-mail: yu-zhou06@mails.tsinghua.edu.cn} &
{\footnotesize E-mail: bzhu@math.tsinghua.edu.cn}
\end{tabular}
}
\bigskip

\today

\end{center}
\def\G{{\Gamma}}
\def\s{\stackrel}
\def\gama{\gamma}
\def\Longrightarrow{{\longrightarrow}}
\def\P{{\cal P}}
\def\A{{\cal A}}
\def\C{\mathcal{C}}
\def\D{\mathcal{D}}
\def\H{\mathcal{H}}
\def\m{\textbf{ M}}
\def\t{{\tau }}
\def\s{\stackrel}
\def\mathbb{\NN}

\def\Hom{\mbox{Hom}}
\def\Ext{\mbox{Ext}}
\def\ind{\mbox{ind}}
\def\deg{\mbox{deg}}

\renewcommand{\mod}{\operatorname{mod}\nolimits}
\newcommand{\add}{\operatorname{add}\nolimits}
\newcommand{\Rad}{\operatorname{Rad}\nolimits}
\newcommand{\RHom}{\operatorname{RHom}\nolimits}
\newcommand{\uHom}{\operatorname{\underline{Hom}}\nolimits}
\newcommand{\End}{\operatorname{End}\nolimits}
\renewcommand{\Im}{\operatorname{Im}\nolimits}
\newcommand{\Ker}{\operatorname{Ker}\nolimits}
\newcommand{\Coker}{\operatorname{Coker}\nolimits}
\renewcommand{\r}{\operatorname{\underline{r}}\nolimits}
\def \text{\mbox}

\hyphenation{ap-pro-xi-ma-tion}

\begin{abstract} We study the cluster combinatorics of
$d-$cluster tilting objects in $d-$cluster categories. Using
mutations of maximal rigid objects in $d-$cluster categories, which
are defined in a similar way to mutations for $d-$cluster tilting
objects, we prove the equivalences between $d-$cluster tilting
objects, maximal rigid objects and complete rigid objects. Using the
chain of $d+1$ triangles of $d-$cluster tilting objects in
\cite{IY}, we prove that any almost complete $d-$cluster tilting
object has exactly $d+1$ complements, compute the extension groups
between these complements, and study the middle terms of these $d+1$
triangles. All results are the extensions of corresponding results
on cluster tilting objects in cluster categories established for
$d-$cluster categories in \cite{BMRRT}. They are applied to the
Fomin-Reading generalized cluster complexes of finite root systems
defined and studied in \cite{FR2} \cite{Th} \cite{BaM1,BaM2}, and to
that of infinite root systems \cite{Zh3}.

\end{abstract}


\medskip


\medskip

\section{Introduction}
Cluster categories are introduced by
Buan-Marsh-Reineke-Reiten-Todorov \cite{BMRRT} for a categorified
understanding of cluster algebras introduced by Fomin-Zelevinsky in
\cite{FZ1,FZ2}, see also \cite{CCS} for type $A_n$. We refer
\cite{FZ3} for a survey on cluster algebras and their combinatorics,
see also \cite{FR1}. Cluster categories are the orbit categories
$\mathcal{D}/\tau^{-1}[1]$ of derived categories of hereditary
categories by the automorphism group $<\tau^{-1}[1]>$ generated by
the automorphism $\tau^{-1}[1]$. They are triangulated categories
\cite{Ke}. Cluster categories, on the one hand, provide a successful
model for acyclic cluster algebras and their cluster combinatoric;
see, for example, \cite{BMRRT}, \cite{BMR}, \cite{CC},
\cite{CK1,CK2}, \cite{IR}, \cite{Zh1,Zh2}; on the other hand,
 they replace module categories as a new generalization of the classical
tilting theory, see, for example, \cite{KR1,KR2}, \cite{IY},
\cite{KZ}. Cluster tilting theory and its combinatorics are the
essential ingredients in the connection between quiver
representations and cluster algebras, and have now become a new part
of tilting theory in the representation theory of algebras; we refer
to the surveys \cite{BM}, \cite{Rin}, \cite{Re} and the references
there for recent developments and background on cluster tilting
theory.
\medskip

Let $H$ be a finite dimensional hereditary algebra over a field $K$
with $n$ non-isomorphic simple modules, and let $\C(H)$ be the
corresponding cluster category.  In a triangulated category, there
are three possible kinds of rigid objects: cluster tilting (maximal
$1-$orthogonal in the sense of Iyama \cite{I}), maximal rigid, and
complete rigid. It is well-known that they are not equivalent to
each other in general \cite{BIKR} \cite{KZ}. But in the cluster
category $\C(H)$,  they are equivalent \cite{BMRRT}. Compared with
 classical tilting modules, cluster tilting objects in cluster
categories have nice properties \cite{BMRRT}. For example, any
almost complete cluster tilting object in a cluster category can be
completed to a cluster tilting object in exactly two ways, but in
mod$H$, there are at most two ways to complete an almost complete
basic tilting module. Moreover, the two complements $M$, $M^{\ast}$
of an almost complete basic cluster tilting object $\bar{T}$ are
connected by two triangles
\begin{displaymath}M^{\ast}{\longrightarrow}B{\longrightarrow}M{\longrightarrow}M^{\ast}[1]\end{displaymath}
\begin{displaymath}M{\longrightarrow}B '{\longrightarrow}M^{\ast}{\longrightarrow}M[1]\end{displaymath}
in $\C(H)$, where respectively, $B{\rightarrow}M$ and $B
'{\rightarrow}M^{\ast}$ are minimal right
add$\bar{T}-$approximations of $M$ and $M^{\ast}$ in $\C(H)$. It
follows that $M$ and $M^{\ast}$ satisfy the condition
 dim$_{D_{M}}\Ext^1_{\C (H)}(M,M^{\ast})\newline =1=\mbox{dim}_{D_{M^{\ast}}}\Ext^1_{\C (H)}(M^{\ast},M)$,
where $D_M$ (or $D_{M^*})$ is the endomorphism division ring of $M$
(resp. $M^*$). Conversely, if two indecomposable rigid objects $M$,
$M^{\ast}$ satisfy the condition above, one can find an almost
complete cluster-tilting object $\bar{T}$ such that $M$ and
$M^{\ast}$ are the two complements of $\bar{T}$. In this case,
$\bar{T}\oplus M^{\ast}$ is called a mutation of $\bar{T}\oplus M$.
Any two cluster-tilting objects are connected through mutations,
 provided that the ground field $K$ is algebraically closed.

\medskip

 Keller \cite{Ke} introduced $d-$cluster categories $\mathcal{D}/\tau^{-1}[d]$ as a
generalization of cluster categories for $d\in \mathbf{N}$. They are
studied recently in \cite{Th}, [Zh3] [BaM1, BaM2], [KR1, KR2],
\cite{IY}, \cite{HoJ1, HoJ2}, \cite{J},  \cite{Pa}, \cite{ABST},
\cite{T}, \cite{Wr}. $d-$cluster categories are triangulated
categories with Calabi-Yau dimension $d+1$ \cite{Ke}. When $d=1$,
 ordinary cluster categories are recovered.

\medskip

The aim of this paper is to study the cluster tilting theory in
$d-$cluster categories. It is motivated by two factors. First, since
some properties of cluster tilting objects in cluster categories do
not hold in general in this generalized setting (for example, the
endomorphism algebras of $d-$cluster tilting objects are not again
Goreistein algebras of dimension at most $d$ in general \cite{KR1}),
 one natural question is to see whether other properties of cluster tilting
  objects hold in $d-$cluster categories. Second,
 in [Zh3] we use $d-$cluster categories to define a generalized
cluster complexes of the root systems of the corresponding Kac-Moddy
Lie algebras (see also \cite{BMRRT} and \cite{Zh1} for a quiver
approach of cluster complexes). When $H$ is of finite representation
type, these complexes are the same as those defined by Fomin-Reading
\cite{FR2} using the combinatorics of the root systems, see also
\cite{Th}. We need the combinatorial properties of $d-$cluster
tilting objects for these generalized cluster complexes.
\medskip

In [Zh3], the second author of this paper proved that any basic
$d-$cluster tilting object in a $d-$cluster category
$\mathcal{C}_d(H)$ contains exactly $n$ indecomposable direct
summands, where $n$ is the number of non-isomorphic simple
$H-$modules, and that the number of complements of an almost
complete $d-$cluster tilting object is at least $d+1$. The present
article is a completion of the result from [Zh3] mentioned above.
Furthermore, it can be viewed as a generalization to $d-$cluster
categories of (almost) all the results for cluster categories in
[BMRRT].
\medskip

The paper is organized as follows: In Section 2, we recall and
collect some notion and basic results needed in this paper. In
Section 3, we prove that the $d-$cluster tilting objects in
$d-$cluster categories are equivalent to the maximal rigid objects,
and also to the complete rigid objects (i.e. rigid objects
containing $n$ non-isomorphic indecomposable direct summands, where
$n$ is the number of simple modules over the associated hereditary
algebra). In the Dynkin case, this equivalence was proved in
\cite{Th} using the fact that every indecomposable object is rigid.
In Section 4, we compare two chains of
 $d+1$ triangles, from \cite{Zh3} and \cite{IY} respectively, in order to
 prove that a basic almost complete $d-$cluster tilting object has
exactly $d+1$ non-isomorphic complements, which are connected by
these $d+1$ triangles. The extension groups between the complements
 of an almost complete $d-$cluster tilting object are computed
explicitly, and a necessary and sufficient condition for $d+1$
indecomposable rigid objects to be the complements of an almost
complete $d-$cluster tilting object is obtained in Section 5. In
Section 6, for an almost complete $d-$cluster tilting object, the
middle terms of the $d+1$ triangles which are connected by the $d+1$
complements are proved to contain no  direct summands common to them
 all. In the final section, we give an application of the results proved
in these previous sections to the generalized cluster complexes
defined by Fomin-Readings \cite{FR2}, studied in \cite{Th}, and
\cite{Zh3}, and show that all the main properties of these
generalized cluster complexes of finite root system in \cite{FR2}
\cite{Th} hold also for the generalized cluster complexes of
arbitrary root systems defined in \cite{Zh3}.

\medskip

After completing and submitting this work, we saw Wralsen's paper
[W] (arXiv 0712.2870). The fact that maximal $d-$rigid objects and
$d-$cluster tilting objects coincide and that almost complete
$d-$cluster tilting objects have $d+1$ complements, have also been
proved independently in [W], with different proofs.

\bigskip

\section{Basics on $d-$ cluster categories}

In this section, we collect some basic definitions and fix notation
that we will use throughout the paper.

Let $H$ be a finite dimensional hereditary algebra over a field $K$.
We denote by $\mathcal{H}$ the category of finite dimensional
modules over $H$. It is a hereditary abelian category \cite{DR}. The
subcategory of $\mathcal{H}$ consisting of isomorphism classes of
indecomposable $H-$modules is denoted by $\ind \mathcal{H}$. The
bounded derived category of $\H$ will be denoted by $D^b(H)$ or
$\D$. We denote the non-isomorphic indecomposable projective
representations in $\mathcal{H}$ by $P_1, \cdots ,P_n$, and the
simple representations with dimension vectors $\alpha _1, \cdots,
\alpha _n$ by $E_1, \cdots, E_n$. We use $D(-)$ to denote
$\mbox{Hom}_K(-,K)$ which is a duality operation in $\mathcal{H}.$

\medskip

The derived category $\cal{D}$ has Auslander-Reiten triangles, and
the Auslander-Reiten translate $\tau$ is an automorphism of
$\mathcal{D}$. Fix a positive integer $d$, and denote by
$F_d=\tau^{-1}[d]$,
 it is an automorphism of $\mathcal{D}$. The $d-$cluster category of $H$ is defined in \cite{Ke};
 we denote by $\mathcal{D}/ F_d$ the corresponding factor category.
 Its objects are by definition the $F_d$-orbits of objects in
$\cal{D}$, and the morphisms are given by
$$\Hom_{\mathcal{D}/F_d}(\widetilde{X},\widetilde{Y}) =
\oplus_{i \in \mathbf{Z}}
 \Hom_{\mathcal{D}}(X,F_d^iY).$$
Here $X$ and $Y$ are objects in $\cal{D}$, and $\widetilde{X}$ and
$\widetilde{Y}$ are the corresponding objects in $\mathcal{D}/F_d$
(although we shall sometimes write such objects simply as $X$ and
$Y$).

\begin{defn}\cite{Ke}\cite{Th}  The orbit category $\mathcal{D}/F_d$ is called the $d-$cluster category of $\mathcal{H}$
(or of $H$), and is denoted by $\mathcal{C}_d(\mathcal{H})$,
 or sometimes by $\mathcal{C}_d(H)$.\end{defn}

By \cite{Ke}, the $d-$cluster category is a triangulated category
with shift functor $[1]$ induced by the shift functor in
$\mathcal{D}$; the projection $\pi: \mathcal{D}\longrightarrow
\mathcal{D}/F$ is a triangle functor. When $d=1$, this orbit
category is called the cluster category of $\mathcal{H}$, and
 denoted by $\mathcal{C}(\mathcal{H})$, or sometimes by
$\mathcal{C}(H)$.

$\mathcal{H}$ is a full subcategory of $\mathcal{D}$ consisting of
complexes concentrated in degree $0$. Passing to
$\mathcal{C}_d(\mathcal{H})$ by the projection $\pi$,  $\mathcal{H}$
is a (possibly not full) subcategory of
$\mathcal{C}_d(\mathcal{H})$, and $\mathcal{C(H)}$ is also a
(possibly not full) subcategory of $\mathcal{C}_d(\mathcal{H})$. For
any $i\in \mathbf{Z}$, we use $(\mathcal{H})[i]$ to denote the copy
of $\mathcal{H}$ under the $i-$th shift $[i]$, considered as a
subcategory of $\mathcal{C}_d(\mathcal{H}).$ Thus,
$(\mbox{ind}\mathcal{H})[i]=\{ M[i]\ | \ M\in \mbox{ind}
\mathcal{H}\ \}$. For any object $M$ in
$\mathcal{C}_d(\mathcal{H})$, let add$M$ denote the full subcategory
of $\mathcal{C}_d(\mathcal{H})$ consisting of direct summands of
direct sums of copies of $M$.

For $X, Y\in \mathcal{C}_d(\mathcal{H})$, we will use Hom$(X,Y)$ to
denote the Hom-space Hom$_{\mathcal{C}_d(\mathcal{H})}(X,Y)$ in the
$d-$cluster category $\mathcal{C}_d(\mathcal{H})$ throughout the
paper. We define Ext$^i(X,Y)$ to be Hom$(X,Y[i]).$
\medskip

We summarize some known facts about $d-$cluster categories
\cite{BMRRT, Ke}, see also \cite{Zh3}.

\begin{prop}\label{pr}

\begin{enumerate}
\item $\mathcal{C}_d(\mathcal{H})$ has Auslander-Reiten triangles
and Serre functor $\Sigma =\tau [1]$, where $\tau$ is the
AR-translate in $\mathcal{C}_d(\mathcal{H})$, induced from
 the AR-translate in $\mathcal{D}$. \item $\mathcal{C}_d(\mathcal{H})$ is
a Calabi-Yau category of CY-dimension $d+1$.
\item $\mathcal{C}_d(\mathcal{H})$ is a Krull-Remak-Schmidt category.
\item $\mathrm{ind}\mathcal{C}_d(\mathcal{H})=\bigcup_{i=0}^{i=d-1}(\mathrm{ind}\mathcal{H})[i]\bigcup\{P_j[d]\ | \ 1\leq j\leq n\}.$
\end{enumerate}
\end{prop}

\begin{proof} See \cite{Zh3}.
\end{proof}

 Using Proposition \ref{pr}, we can define the degree for every
indecomposable object in $\mathcal{C}_d(\mathcal{H})$ as follows
\cite{Zh3}:

\begin{defn}\label{degree} For any indecomposable object $X\in \mathcal{C}_d(\mathcal{H})$, we call the non-negative integer
$\mathrm{min }\{ k\in \mathbf{Z}_{\geq 0} \  | \  X \cong M[k]
\mbox{ in } \mathcal{C}_d(\mathcal{H}), \mbox{ for  some  } M\in
\mathrm{ind}\mathcal{H}\ \}$ the degree of $X$, denoted by $\deg X$.
If $\deg X=k,k=0, \cdots, d-1$, we say that $X$ is of color $k+1$;
if $\deg X=d$, we say that $X$ is of color $1$.
\end{defn}

By Proposition \ref{pr}, any indecomposable object $X$ of degree $k$
is isomorphic to $M[k]$ in $\mathcal{C}_d(\mathcal{H})$, where $M$
is an indecomposable representation in $\mathcal{H}$, $0\leq \deg
X\leq d$, $X$ has degree $d$ if and only if $X\cong P[d]$ in
$\mathcal{C}_d(\mathcal{H})$ for some indecomposable projective
object $P\in \mathcal{H}$, and $X$ has degree $0$ if and only if
$X\cong M[0]$ in $\mathcal{C}_d(\mathcal{H})$ for some
indecomposable object $M\in \mathcal{H}$. Here $M[0]$ denotes
 the object $M$ of $\mathcal{H}$, considered as a complex concentrated
in degree $0$.

Now we recall the notion of $d-$cluster tilting objects from
\cite{KR1}, \cite{Th}, \cite{Zh3}, \cite{IY}. This notion is
equivalent to the "maximal $d-$orthogonal subcategories" of Iyama
\cite{I, IY}.

\medskip

\begin{defn} \label{defmaxorth}
Let $\mathcal{C}_d(\mathcal{H})$ be the d-cluster category.

\begin{enumerate}
\item An object $X$ in $\C_{d}(\H)$ is called rigid if
$\Ext^i(X,X)=0$, for all $1\leq i\leq d$.

\item An object $X$ in $\C_{d}(\H)$ is called maximal rigid if it
satisfies the property: $Y \in \add X$ if and only if
$\Ext^i(X\bigoplus Y,X\bigoplus Y)=0$ for all $1\leq i\leq d$.

\item An object $X$ in $\C_{d}(\H)$ is called completely rigid if it
contains exactly $n$ non-isomorphic indecomposable direct summands.

\item An object $X$ in $\C_{d}(\H)$ is called d-cluster tilting if
it satisfies the property that $Y \in \add X$ if and only if
$\Ext^i(X,Y)=0$ for all $1\leq i\leq d$.

\item An object $X$ in $\C_{d}(\H)$ is called an almost complete
d-cluster tilting if there is an indecomposable object $Y$ with $Y
\notin \add X$ such that $X \bigoplus Y$ is a d-cluster tilting
object. Such $Y$ is called a complement of the almost complete
$d-$cluster tilting object.
\end{enumerate}

\end{defn}

For a basic d-cluster tilting object $T$ in $\C_d(H)$, an
indecomposable object $X_0 \in \mbox{add}T$ and its complement $X$
such that $X_0 \bigoplus X=T$, then there is a triangle in
$\C_d(H)$:
\begin{displaymath}X_1\s{g}{\longrightarrow}B_0\s{f}{\longrightarrow}X_0{\longrightarrow}X_1[1],\end{displaymath}
where $f$ is the minimal right add$X-$approximation of $X_0$ and $g$
is the minimal left add$X-$approximation of $X_1$. It is easy to see
that $T':=X_1 \bigoplus X$ is a basic d-cluster tilting object
(compare \cite{IY}). We call $T'$ is a mutation of $T$ in the
 direction of $X_0$. We call two $d-$cluster tilting objects $T, \
T'$ mutation equivalent provided that there are finitely many
$d-$cluster tilting objects $T_1(=T),T_2,\cdots, T_n(=T')$ such that
$T_{i+1}$ is a mutation of $T_i$ for any $1\le i\le n-1$.
\medskip

 From the proof of Theorem 4.6 in \cite{Zh3}, we know that every
$d-$cluster tilting object is mutation equivalent to a $d-$cluster
tilting object in $\mathcal{H}[0]$.

\medskip

The following results are proved in [Zh3].

\begin{prop}

\begin{enumerate}
\item Any indecomposable rigid object $X$ in
$\mathcal{C}_d(\mathcal{H})$ is either of the form $M[i]$, where
 $M$ is a rigid  module (i.e. Ext$^1_H(M,M)=0$) in
$\mathcal{H}$ and $0\leq i\leq d-1$, or of the form $P_j[d]$ for
some $1\leq j\leq n$. In particular, if $\G$ is a Dynkin graph, then
any indecomposable object in $\mathcal{C}_d(\mathcal{H})$ is rigid.
\item Suppose $d \geq 2$. Then $\End_{\C_d(H)}(X)$ is a division
algebra for any indecomposable rigid object $X$.

\item Let $d \geq 2$ and X=M$[i]$, Y=N$[j]$ be indecomposable
objects of degree i,j respectively in $\C_d(H)$. Suppose that
$\Hom(X,Y) \neq 0$. Then one of the following holds:

$(1)$We have $i=j$ or $j-1$ $($provided $j \geq 1)$;

$(2)$We have $i=0$, $i=d$ $($and $M=P)$ or $d-1$ $($provided
$j=0)$.

\item Let $d \geq 2$ and $M,N \in \H$. Then any non-split triangle
between $M[0]$ and $N[0]$ in $\C_d(H)$ is induced from a non-split
exact sequence between M and N in $\H$.
\end{enumerate}

\end{prop}

\section{Equivalence of $d-$cluster tilting objects and maximal rigid objects}

The equivalence between cluster tilting objects and maximal rigid
objects in cluster categories was proved in \cite{BMRRT}.  For
$d-$cluster categories, in the simply laced Dynkin case, the
equivalence of $d-$cluster tilting objects and maximal rigid objects
is easily obtained because any indecomposable object is rigid
(compare \cite{Th}). We will now prove it for arbitrary $d-$cluster
categories. From the proof of Theorem 4.6 in \cite{Zh3}, we know
that every $d-$ cluster tilting object is mutation equivalent to one
in $\mathcal{H}[0]$. If there is a similar result for mutations of
 maximal rigid objects, then
we can get the equivalence by the obvious equivalence between
$d-$cluster tilting objects and maximal rigid objects in
$\mathcal{H}[0]$ (both are tilting modules in mod$H$).

\begin{lem}\label{maxmut}

Let $d \geq 2$, $T=X\bigoplus X_0$ be a basic maximal rigid object
in $\mathcal{C}_d(\mathcal{H})$ and $X_0$ an indecomposable object.
Then there are $d+1$ triangles
\begin{displaymath}(*) \ \ \  X_{i+1}\s{g_i}{\longrightarrow}T_i\s{f_i}{\longrightarrow}
X_i\s{\delta_i}{\longrightarrow}X_{i+1}[1],\end{displaymath}

where $T_i \in \add X$, $f_i$ is the minimal right
add$X-$approximation of $X_i$, $g_i$ is the minimal left
add$X-$approximation of $X_{i+1}$, all the $X\bigoplus_{}X_i$ are
maximal rigid objects, and all $X_i$ are distinct up to isomorphisms
for $i=0, \cdots, d$.

\end{lem}

\begin{proof}  First we prove that there is a triangle
\begin{displaymath}  X_{1}\s{g_0}{\longrightarrow}
T_0\s{f_0}{\longrightarrow}X_0\s{\delta_0}{\longrightarrow}X_1[1],\end{displaymath}

where $T_0 \in \mbox{add}X$, $f_0$ is the minimal right
$\mbox{add}X-$approximation of $X_0$, $g$ is the minimal left
$\mbox{add}X-$approximation of $X_{1}$,  and $X\bigoplus_{}X_1$ is a
maximal rigid object.

Let $T_0\s{f_0}{\longrightarrow}X_0$ be the minimal right
$\mbox{add}X-$approximation of $X_0$, and let
\begin{displaymath}(1) \ \ \ X_1\s{g_0}{\longrightarrow}T_0\s{f_0}{\longrightarrow}X_0\s{\delta_0}{\longrightarrow}X_1[1]\end{displaymath}
be the triangle into which $f$ embeds. By the discussion in
\cite{BMRRT}, one can easily check that $g_0$ is the minimal left
$\mbox{add}X-$approximation of $X_1$, $X_1$ is indecomposable and
$X_1\notin \add X$. By applying $\Hom(X,-)$ to the triangle, we have
$\Ext^{i}(X,X_1)=0$, for $1\leq i\leq d$ (for $i=1$, because $f$ is
the minimal right add$X$-approximation of $X_0$). By applying
$\Hom(X_0,-)$ to the triangle, we get $\Ext^{i}(X_0,X_0)\cong
\Ext^{i+1}(X_0,X_1)$, for $1\leq i\leq d-1$. By applying
$\Hom(-,X_1)$ to the triangle, we have $\Ext^{i}(X_1,X_1) \cong
\Ext^{i+1}(X_0,X_1)$, for $1\leq i \leq d-1$. So
 $\Ext^{i}(X_1,X_1)\cong \Ext^{i}(X_0,X_0)=0$ for $1 \leq i \leq
d-1$. Since $\C_d(\mathcal{H})$ is a Calabi-Yau category of
CY-dimension d+1, $\Ext^{d}(X_1,X_1)\cong D\Ext^{1}(X_1,X_1)=0$. We
claim that $X\bigoplus X_1$ is a maximal rigid object. If not, we
have an indecomposable object $Y_1 \notin \mbox{add}(X\bigoplus
X_1)$, such that $X\bigoplus X_1\bigoplus Y_1$ is a rigid object.
Then we have a triangle
\begin{displaymath}(2) \ \ \ Y_1\s{\psi}{\longrightarrow}T_1\s{\varphi}{\longrightarrow}Y_0{\longrightarrow}X_1[1],\end{displaymath}
where $\psi$ is the minimal left $\mbox{add}X-$approximation of
$Y_1$. It is easy to prove that $\varphi$ is the minimal right
$\mbox{add}X-$approximation of $Y_0$, $Y_0\notin \mbox{add}X$, and
$\Ext^{i}(Y_0,X\bigoplus Y_0)=0$ for $1\leq i\leq d$. We will prove
that $\Ext^{i}(Y_0,X_0)=0$ for $1\leq i\leq d$; then $Y_0 \cong X_0$
 due to the fact that $X\oplus X_0$ is a maximal rigid object.
 By applying $\Hom(-,Y_1)$ to the first triangle, we
have $0=\Ext^{i}(X_1,Y_1)\cong \Ext^{i+1}(X_0,Y_1)$ for $1\leq i\leq
d-1$. By applying $\Hom(X_0,-)$ to the second triangle, we have
$\Ext^{i}(X_0,Y_0)\cong \Ext^{i+1}(X_0,Y_1)=0$ for $1\leq i\leq
d-1$. So we have $\Ext^{i}(X_0,Y_0)=0$ for $1\leq i\leq d-1$, and
 thus $\Ext^{i}(Y_0,X_0)=0$ for $2\leq i\leq d$. By applying
$\Hom(-,X_1)$ to the second triangle, we have
$0=\Ext^{1}(Y_1,X_1)\cong \Ext^{2}(Y_0,X_1)$. By applying
$\Hom(Y_0,-)$ to the first triangle, we have $\Ext^{1}(Y_0,X_0)\cong
\Ext^{2}(Y_0,X_1)=0$. So $\Ext^{1}(Y_0,X_0)=0$. In all,
$\Ext^{i}(Y_0,X_0)=0$ for $1 \leq i \leq d$. Therefore $Y_0\cong
X_0$ which induces an isomorphism between the triangles (1) and (2).
Then $Y_1\cong X_1$, a contradiction. This proves that $X\oplus X_1$
is a maximal rigid object.

Second we repeat this process to get $d+1$ triangles
\begin{displaymath}(*) \ \ \ X_{i+1}\s{g_i}{\longrightarrow}T_i
\s{f_i}{\longrightarrow}X_i\s{\delta_i}{\longrightarrow}X_{i+1}[1],\end{displaymath}
where $T_i \in \add X$, $f_i$ is the minimal right
$\mbox{add}X-$approximation of $X_i$, $g_i$ is the minimal left
$\mbox{add}X-$approximation of $X_{i+1}$, and all the
 $X\bigoplus_{}X_i$ are maximal rigid objects.

Third it is easy to see that $\delta_d[d]\delta_{d-1}[d-1]\cdots
\delta_1[1]\delta_0\not=0$ (similar as that in Corollary 4.5 in
\cite{Zh3}). In particular, $\mathrm{Hom}(X_i,X_j[j-i])\not=0$ and
$X_i\ncong X_j, \forall 0\leq i<j\leq d.$ This finishes the proof.
\end{proof}

With the help of Lemma \ref{maxmut}, one can define mutations of
maximal rigid objects similar to those of $d-$cluster tilting
objects: Let
\begin{displaymath} \ \ \ X_{i+1}\s{g_i}{\longrightarrow}T_i\s{f_i}{\longrightarrow}
X_i\s{\delta_i}{\longrightarrow}X_{i+1}[1]\end{displaymath} be the
$i-$th triangle in Lemma \ref{maxmut}. We say that each of the
maximal rigid objects $X\oplus X_i$, for $i=1, \cdots, d$, is a
mutation of the maximal rigid object $X\oplus X_0$. A maximal rigid
object $T$ is mutation equivalent to a maximal rigid object $T'$
provided that there are finitely many maximal rigid objects
$T_1(=T), T_2,\cdots, T_{n-1},T_n(=T')$ such that $T_i$ is a
mutation of $T_{i-1}$ for any $i$.

\medskip

\begin{lem}\label{mutto0}
Let $d \geq 2$, $T=X\bigoplus X_0$ be a maximal rigid object and
$X_0$ be an indecomposable object. Then $T$ is mutation equivalent
to a maximal rigid object in $\mathcal{H}[0]$.
\end{lem}

\begin{proof} In the proof of Theorem 4.6 in \cite{Zh3}, we proved that
any $d-$cluster tilting object is mutation equivalent to a
$d-$cluster tilting object in $\mathcal{H}[0]$. The same proof works
here (with the help of Lemma 3.1), after replacing $d-$cluster
tilting objects by maximal rigid objects. We omit the details and
refer to the proof of Theorem 4.6. in \cite{Zh3}.
\end{proof}

Now we prove the main result in this section.

\begin{thm} \label{charmonoepi}
Let $X$ be a basic rigid object in the d-cluster category
$\C_d(\mathcal{H})$. Then the following statements are equivalent:

\begin{enumerate}

\item $X$ is a d-cluster tilting object.

\item $X$ is a maximal rigid object.

\item $X$ is a complete rigid object, i.e. it contains exactly n indecomposable summands.

\end{enumerate}

\end{thm}

\begin{proof} We suppose that $d > 1$; the same statement was proved for $d = 1$ in \cite{BMRRT}.
We prove that the first two conditions are equivalent. A d-cluster
tilting object must be a maximal rigid object by definition. Now we
assume $X$ is a maximal rigid object. Then $X$ is mutation
equivalent to a maximal rigid object $T'[0]$ in $\mathcal{H}[0]$ by
Lemma \ref{mutto0}. We have that $\Ext^k(T'[0],T'[0])\cong
\Ext_{\mathcal{D}}^k(T'[0],T'[0])\cong \Ext_{\mathcal{H}}^k(T',T')$,
$k=1, \cdots, d-1$ and $\Ext^d(T'[0],T'[0])\cong
D\Ext(T'[0],T'[0])\cong D\Ext_{\mathcal{H}}(T',T')$. So $T'$ is a
maximal rigid module in $\mathcal{H}$. Hence $T'$ is a tilting
module, and thus $T'[0]$ is a $d-$cluster tilting object. Therefore
$T$ is a $d-$cluster tilting object, since it is mutation equivalent
to the $d-$cluster tilting object $T'[0]$.

Now we prove that the last two conditions are equivalent. In
\cite{Zh3}, we know that every basic d-cluster tilting object has
exactly $n$ indecomposable summands. Conversely, any basic rigid
object with $n$ indecomposable summands will be a basic maximal
rigid object, since otherwise it can be extended to a basic maximal
rigid object that contains at least $n+1$ indecomposable summands.
This is a contradiction.
\end{proof}

This theorem immediately yields the following important conclusion.

\begin{cor}
Let $X$ be a rigid object in $\C_{d}(\H)$. Then there exists an
object $Y$ such that $X \bigoplus Y$ is a d-cluster tilting object.
\end{cor}

\section{Complements of almost complete basic $d$-cluster tilting objects}

The number of complements of an almost complete cluster tilting
object in a cluster category $\C(\mathcal{H})$ is exactly two
\cite{BMRRT}. From Corollary 4.5 in [Zh3], we know that the number
of complements of an almost complete $d-$cluster tilting object is
at least $d+1$. In this section, we will prove it is exactly $d+1$.
\medskip

Let $T=X\bigoplus X_0$ be a basic $d-$cluster tilting object in
$\mathcal{C}_d(\mathcal{H})$, and $X$ an almost complete $d-$cluster
tilting object. By Theorem 4.4 in \cite{Zh3} and Theorem 3.10 in
\cite{IY}, we have the following two chains of $d+1$ triangles:

\begin{displaymath}(*)\ \  X_{i+1}\s{g_i}{\longrightarrow}B_i\s{f_i}{\longrightarrow}
X_i\s{\delta_i}{\longrightarrow}X_{i+1}[1],\end{displaymath} where
 for $i=0, 1, \cdots, d$, $B_i \in \mbox{add}X$, the map $f_i$ is the minimal right add$X-$approximation of
$X_i$ and $g_i$ is the minimal left add$X-$approximation of
$X_{i+1}$.\vspace{10pt}

\begin{displaymath}(**)\  \ X'_{i+1}\s{b_i}{\longrightarrow}C_i
\s{a_i}{\longrightarrow}X'_i\s{c_i}{\longrightarrow}X'_{i+1}[1],\end{displaymath}
where for $i=0, 1, \cdots, d$, $C_i \in \mbox{add}T$, the map $a_i$
is the minimal right add$T-$approximation of $X'_i$ (except $a_0$,
 which is the sink map of $X'_0$ in add$T$) and $b_i$ is the minimal left
add$T-$approximation of $X'_{i+1}$ (except $b_d$, which is the
source map of $X'_d$ in add$T$), and $X'_0=X'_{d+1}=X_0$.
\medskip

In \cite{IY}, the authors show that $X_0 \notin \add(\bigoplus_{0
\leq i \leq d}C_i)$ is a sufficient condition for an almost complete
$d-$cluster tilting object to have exactly $d+1$ complements. The
main aim of this section is to prove that $B_i=C_i$ for all $0 \leq
i \leq d$, which implies this sufficient condition. We will first
study the properties of the degree of an indecomposable object in
$\C_d(\mathcal{H})$ which is a useful tool for studying rigid
objects in $d-$cluster categories.

\begin{lem}

Let $X_i, \ 0\le i\le d ,$  be the objects appearing in the
triangles in $(*)$.  If $\deg X_0=0$, then

(1) $\deg X_1=0$, $d$ or $d-1$, and

(2) $\deg X_i \geq d-i$, for any $2 \leq i \leq d$.

\end{lem}

\begin{proof}

(1) We have the fact that $\Hom(X_0,X_1[1])=\Ext(X_0,X_1) \neq 0$.
If $0< \deg X_1<d-1$ (which implies $d \geq 3$), then $2 \leq \deg
X_1[1] \leq d-1$ and $\Hom(X_0,X_1[1])=0$ by Proposition 2.5(3).
This is a contradiction.

(2) If $\deg X_1=0$, then $\deg X_2=d$ or $d-1$ or $d-2$ (because
$X_0$, $X_1$, $X_2$ cannot have the same degree by the proof of
Theorem 4.6 in [Zh3]). Now we prove the assertion that $\deg
X_{i+1}\ge d-(i+1)$ provided that $\deg X_i \geq d-i$ for some $i$
 ($1 \leq i \leq d-1$).
 If $\deg X_{i+1}<d-(i+1)$, then $1 \leq \deg X_{i+1}[1]<d-i$, which
implies $d \geq 2$, and then $\Hom(X_i,X_{i+1}[1])=0$ by Proposition
2.5. This contradicts the fact $\Ext(X_i,X_{i+1}) \neq 0$. So by
induction on $i$, we get the statement (2).

\end{proof}

\begin{lem}\label{computer}

Let $d \geq 2$ and $X=M[i]$, $Y=N[j]$ be indecomposable objects of
degree i, j respectively in $\C_{d}(\H)$. Suppose that $0 \leq j+k-i
\leq d-1$. Then

(1) $\Hom(X,Y[k]) \cong \Hom_{\D}(X,Y[k])$, and

(2) $\Hom(X,{\t}^{-1}Y[k]) \cong \Hom_{\D}(X,{\t}^{-1}Y[k])$.

\end{lem}

\begin{proof}

(1) $\Hom(X,Y[k])=\bigoplus_{l \in \mathbf{Z}}\Hom_{D}(X,
{\t}^{-l}Y[k+ld])$.

When $l \geq1 $, $\Hom_{\mathcal{D}}(X,{\t}^{-l}Y[k+ld])\cong
\Hom_{\mathcal{D}}({\t}^{l}M,N[k+ld-i+j])=0$, since $k+ld-i+j \geq
ld \geq 2$.

When $l \leq -1$, $\Hom_{\D}(X,{\t}^{-l}Y[k+ld])\cong
D\Hom_{\D}({\t}^{-l-1}N,M[-k-ld+i-j+1])=0$, since $-l-1 \geq 0$ and
$-k-ld+i-j+1 \geq 2-(l+1)d \geq 2$.

It follows that $\Hom(X,Y[k])\cong \Hom_{\D}(X,Y[k])$.

(2) $\Hom(X,{\t}^{-1}Y[k])=\bigoplus_{l \in \mathbf{Z}}\Hom_{D}(X,
{\t}^{-l-1}Y[k+ld])$.

When $l \geq1 $, $\Hom_{\mathcal{D}}(X,{\t}^{-l-1}Y[k+ld])\cong
\Hom_{\mathcal{D}}({\t}^{l+1}M,N[k+ld-i+j])=0$, since $l+1 \geq 2$
and $k+ld-i+j \geq ld \geq 2$.

When $l=-1$, $\Hom_{D}(X,
{\t}^{-l-1}Y[k+ld])=\Hom_{D}(M,N[k-d-i+j])=0$, since $k-d-i+j \leq
-1$.

When $l \leq -2$, $\Hom_{\D}(X,{\t}^{-l-1}Y[k+ld])\cong
D\Hom_{\D}({\t}^{-l-2}N,M[-k-ld+i-j+1])=0$, since $-l-2 \geq 0$ and
$-k-ld+i-j+1 \geq 2-(l+1)d \geq 2$.

It follows that $\Hom(X,{\t}^{-1}Y[k]) \cong
\Hom_{\D}(X,{\t}^{-1}Y[k])$.
\end{proof}

For convenience, we add a triangle below to the triangle chains
$(*)$:
\begin{displaymath}X_0\s{g_{-1}}{\longrightarrow}B_{-1}
\s{f_{-1}}{\longrightarrow}X_{-1}\s{\delta_{-1}}{\longrightarrow}X_{0}[1],\end{displaymath}
where $f_{-1}$ is the right add$X-$approximation and $g_{-1}$ is
 the left add$X-$approximation. Now we prove the main theorem in this
section.

\begin{thm}
Let $d \geq 2$, $T=X\bigoplus X_0$ be a basic $d-$cluster tilting
object in $\mathcal{C}_d(\mathcal{H})$, and $X$ an almost complete
$d-$cluster tilting object. Then there are exactly d+1 complements
$\{X_i\}_{0 \leq i \leq d}$ of $X$, which are connected by the $d+1$
triangles $(*)$. \label{mainthm}
\end{thm}

\begin{proof} The main step in the proof is to show that $X_0 \notin \add C_i$ for $0 \leq
i \leq d$.

For $i=0$ or $i=d$, since $f_0$ is the minimal right
add$X-$approximation of $X_0$ and End$X_0$ is a division ring, for
any map $h \in \Hom(T',X_0)$ that is not a retraction, where $T'$ is
some object in add$T$, there exists $h' \in \Hom(T',B_0)$ such that
$h=f_0h'$. Therefore, $f_0$ is a sink map in add$T$. By the
uniqueness of the sink map, we get $C_0 \cong B_0$, $X_1 \cong X'_1$
and, dually $C_d \cong B_{-1}$, $X_{-1} \cong X'_d$. So $X_0 \notin
\mbox{add}C_0$ and $X_d \notin \mbox{add}C_d$.

For $1 \leq i \leq d-2$ (this implies $d \geq 3$), if $i=1$, by
applying $\Hom(X_0,-)$ to the triangle
$X_2{\longrightarrow}B_1{\longrightarrow}X_1{\longrightarrow}X_2[1]$,
we have the exact sequence
\begin{displaymath}\Hom(X_0,B_1){\longrightarrow}\Hom(X_0,X_1){\longrightarrow}\Ext(X_0,X_2){\longrightarrow}0.\end{displaymath}
We need to prove $\Ext(X_0,X_2)=0$. If not, i.e. $\Ext(X_0,X_2) \neq
0$, then $\Hom(X_0,X_1) \neq 0$. Similarly, by applying
$\Hom(-,X_2)$ to the triangle
$X_1{\longrightarrow}B_0{\longrightarrow}X_0{\longrightarrow}X_1[1]$,
we have the exact sequence
\begin{displaymath}\Hom(X_1,X_2){\longrightarrow}\Ext(X_0,X_2){\longrightarrow}0,\end{displaymath}
so $\Ext(X_0,X_2) \neq 0$ implies $\Hom(X_1,X_2) \neq 0$. We know
that $\Ext(X_0,X_1) \neq 0$ and $\Ext(X_1,X_2) \neq 0$. We may
 assume that the degree of $X_0$ is $0$; then $\deg X_1=0$, $d$ or $d-1$ by
Lemma 4.1. But $\Hom(X_0,X_1) \neq 0$ implies that the degree of
$X_1$ is not $d$ or $d-1$, so it is $0$. For the same reason, $\deg
X_2=0$, which contradicts the fact that $X_0$, $X_1$, and $X_2$
 do not all have the same degree (refer to the proof of Theorem 4.6 in
[Zh3]).

If $2 \leq i \leq d-2$, then by applying $\Hom(X_0,-)$ to the
triangle
$X_{i+1}{\longrightarrow}B_i{\longrightarrow}X_i{\longrightarrow}X_{i+1}[1]$,
we get the exact sequence
\begin{displaymath}\Hom(X_0,B_i){\longrightarrow}\Hom(X_0,X_i){\longrightarrow}\Ext(X_0,X_{i+1}){\longrightarrow}0.\end{displaymath}
We want to prove that $\Hom(X_0,X_i)=0$, which implies
$\Ext(X_0,X_{i+1})=0$. We also assume that the degree of $X_0$ is 0.
 Since $\deg X_i \geq d-i \geq 2$ by Lemma 4.1, it follows that $\Hom(X_0,X_i)=0$.
 So $\Ext(X_0,X_{i+1})=0$, and it follows that $f_i$ is the minimal
right add$T-$approximation of $X_i$. By the uniqueness of the
 minimal approximation map, since $X_1 \cong X'_1$, we get $C_i \cong
B_i$ and $X_{i+1} \cong X'_{i+1}$ for $1 \leq i \leq d-2$, so $X_0
\notin \mbox{add}(\bigoplus_{1 \leq i \leq d-2}C_i)$.

For $i=d-1 \geq 1$ (which implies $d \geq 2$),  we claim that in the
triangle
 $X_d\s{g_{d-1}}{\longrightarrow}B_{d-1}\s{f_{d-1}}{\longrightarrow}X_{d-1}{\longrightarrow}X_d[1]$,
 the morphism $f_{d-1}$ is the minimal right add$(X \bigoplus
X_0)-$approximation of $X_{d-1}$, which is equivalent to the fact
that $\Ext(X_0,X_d)=0$. Suppose that $\deg X_0$=0 and $\deg X_1 \neq
0$ (if $\deg X_0 = \deg X_1=0$, then $\deg X_2 \neq 0$, and we can
 replace $X_0$ by $X_1$). From Lemma 4.1 (2), $\deg X_{d-1} \geq 1$.
 If $\deg X_{d-1} =1$, then $\deg X_d = 1$ or $0$ since
 $\Hom(X_{d-1},X_{d}[1]) \not= 0$. So we divide the calculation of
$\Ext(X_0,X_d)$ into three cases:

\begin{enumerate}

\item The case $\deg X_{d-1} \geq 2$. Then by Proposition 2.1(3) $\Hom(X_0,X_{d-1})=0$, which implies $\Ext(X_0,X_d)=0$.

\item The case $\deg X_{d-1}=1$ and $\deg X_d=1$. By applying $\Hom(X_0,-)$
to the triangle
$X_d{\longrightarrow}B_{d-1}{\longrightarrow}X_{d-1}\s{\delta_{d-1}}{\longrightarrow}X_d[1]$
we get the exact sequence
\begin{displaymath}\Hom(X_0,X_{d-1})\s{\delta_{d-1}^{\ast}}{\longrightarrow}
\Hom(X_0,X_d[1]){\longrightarrow}0,\end{displaymath} where
$\delta_{d-1} \in \Hom(X_{d-1},X_d[1])\cong
\Hom_{\D}(X_{d-1},X_d[1])$ by Lemma 3.2. For any $\varphi \in
\Hom(X_0,X_{d-1})\cong \Hom_{\D}(X_0,X_{d-1})$ by Lemma 4.2, we have
$\delta_{d-1}^{\ast}(\varphi)=\delta_{d-1} \varphi \in
\Hom_{\D}(X_0,X_d[1])=0$. So $\delta_{d-1}^{\ast}=0$. Thus
 $\Ext(X_0,X_d)=0$.

\item The case $\deg X_{d-1}=1$ and $\deg X_d=0$. Consider the triangle
 $X'_d$ ${\longrightarrow}$ $C_{d-1}$ ${\longrightarrow}$$X'_{d-1}$
${\longrightarrow}X'_d[1]$. Since $X_{-1} \cong X_d'$ and $X_{d-1}
\cong X_{d-1}'$, the triangle is $X_{-1}$ ${\longrightarrow}$
$C_{d-1}$ ${\longrightarrow}$ $X_{d-1}$ ${\longrightarrow}$
$X_{-1}[1]$, where $C_{d-1} \in \mbox{add}(X \bigoplus X_0)$.
Analogously, we get a triangle
\begin{displaymath}X_0{\longrightarrow}Y{\longrightarrow}X_d{\longrightarrow}X_0[1],\end{displaymath}
 where $Y \in \mbox{add}(X \bigoplus X_1)$. Since $\deg X_0 = \deg X_d=0$,
 then the degree of the indecomposable summands of $Y$ is zero. But $\deg
X_1 \not= 0$, so $X_1 \notin Y$, that is, $Y \in \mbox{add}X$. By
applying $\Hom(X_0,-)$ to the triangle above, we get the exact
sequence
\begin{displaymath}\Ext(X_0,Y){\longrightarrow}\Ext(X_0,X_d){\longrightarrow}\Ext^2(X_0,X_0){\longrightarrow}X_0[1],\end{displaymath}
 so $\Ext(X_0,X_d)=0$ since $X_0
\bigoplus X$ is a $d-$cluster tilting object.

\end{enumerate}

Then $C_{d-1} \cong B_{d-1}$ so $X_0 \notin \mbox{add}C_{d-1}$.

In all, $X_0 \notin \mbox{add}(\bigoplus_{0 \leq i \leq d}C_i)$,
which satisfies the condition of Corollary 5.9 in \cite{IY}.
Therefore, $X$ has exactly $d+1$ complements in
$\mathcal{C}_d(\mathcal{H})$.

\end{proof}

As a consequence of the proof of the theorem above, we have
\begin{cor}

The corresponding triangles in the chains (*) and (**) are isomorphic.

\end{cor}

Let $d \geq 2$. For a (basic) $d-$cluster tilting object
$T=X\bigoplus X_0$ in $\mathcal{C}_d(\mathcal{H})$ with an almost
complete $d-$cluster tilting object $X$, and for any $i$ between $0$
and $d$, the triangle
\begin{displaymath}
X_{i+1}\s{g_i}{\longrightarrow}B_i\s{f_i}{\longrightarrow}
X_i\s{\delta_i}{\longrightarrow}X_{i+1}[1]\end{displaymath} in $(*)$
 is called the $i-$th connecting triangle of the complements of $X$ with
respect to $X_0$. These $d+1$ triangles form a
$d+1-$Auslander-Reiten triangle starting at $X_0$ (see \cite{IY}).

Similar to the cluster categories in \cite{BMRRT}, one can associate
to $\C _d(\mathcal{H})$ a mutation graph of $d-$cluster tilting
objects: the vertices are the basic $d-$cluster tilting objects, and
 there is an edge between two vertices if the corresponding two basic $d-$cluster
 tilting objects in $\C _d(\mathcal{H})$
have all but one indecomposable summand in common. Exactly as
 in \cite{BMRRT}, we obtain the
 conclusion below, which means that over an algebraically closed field,
 any two d-cluster tilting objects in $\C_d(\mathcal{H})$ can
be connected by a series of mutations.

\begin{prop}
Let $K$ be an algebraically closed field. Given an indecomposable
hereditary $k$-algebra $H$, the associated mutation graph of
$d-$cluster tilting objects in $\C_d(\mathcal{H})$ is connected.
\end{prop}

\section{Relations of complements}

Let $T=X\bigoplus X_0$ be a basic $d-$cluster tilting object in
$\mathcal{C}_d(\mathcal{H})$. The almost complete $d-$cluster object
$X$ has exactly $d+1$ complements $X_i,\ 0\le i\le d$, as shown in
Theorem 4.3. When $d=1$, the extension groups of between $X_0$ and
$X_1$ were computed in \cite{BMRRT}. In this section we will compute
$\Ext^k(X_i,X_j)$. Throughout this section, we assume $d \geq 2$,
and $X$ is a basic almost complete $d-$cluster tilting object, the
$d+1$ complements $X_0,\cdots , X_d$ of $X$ are connected by the
 $d+1$ triangles in $(*)$ in Section 4:

\begin{displaymath}(*)\ \  X_{i+1}\s{g_i}{\longrightarrow}B_i\s{f_i}{\longrightarrow}
X_i\s{\delta_i}{\longrightarrow}X_{i+1}[1],\end{displaymath} where
 for $i=0, 1, \cdots, d$, $B_i \in \mbox{add}X$, $f_i$ is the minimal right
add$X-$approximation of $X_i$ and $g_i$ is the minimal left
 add$X-$approximation of $X_{i+1}$. \vspace{10pt}

\begin{lem}

$\Ext^{i}(X_0,X_i) \cong \Ext(X_0,X_1) \cong \End_{\H}(X_0)$, and
$\Ext^k(X_0,X_i)=0$ for $1 \leq i \leq d, $ and $k\in \{ 1, \cdots
,d\}\backslash \{i\}.$

\end{lem}

\begin{proof}
By applying  $\Hom(X_0,-)$ to the triangles $(*)$ we get the long
exact sequences
\begin{displaymath}\Ext^k(X_0,B_i){\longrightarrow}
\Ext^k(X_0,X_i){\longrightarrow}\Ext^{k+1}(X_0,X_{i+1}){\longrightarrow}\Ext^{k+1}(X_0,B_i),\end{displaymath}
where $i=0, 1, \cdots, d$,  and $k=1, 2, \cdots, d-1$. Since
$\Ext^k(X_0,B_i)=0$ for $0 \leq i \leq d$ and $1 \leq k \leq d$, we
 have $\Ext^k(X_0,X_i)\cong \Ext^{k+1}(X_0,X_{i+1})$ for $0 \leq i \leq
d$ and $1 \leq k \leq d-1$. So $\Ext^{i+1}(X_0,X_{i+1})\cong
\Ext^{i}(X_0,X_{i})$, for $1 \leq i \leq d-1$. Hence we get the left
equation by induction on $i$. Applying $\Hom(X_0,-)$ to the triangle
$X_1{\longrightarrow}B_0{\longrightarrow}X_0\s{\delta_0}{\longrightarrow}X_1[1]$
 induces the exact sequence
\begin{displaymath}\Hom(X_0,X_0)\s{\delta_0^{\ast}}{\longrightarrow}\Ext(X_0,X_1){\longrightarrow}0.\end{displaymath}
Since $\Hom(X_0,X_0)$ is a division algebra for $d \geq 2$, it
 follows that ${\delta_0^{\ast}}(\varphi)=\delta_0 \varphi $ is non-zero
for any non-zero map $\varphi$ in End$X_0$, which must therefore be
an isomorphism of $X_0$. Then ${\delta_0^{\ast}}$ is a monomorphism
and hence an isomorphism. This gives the first part of the lemma.

For the second part, if $i < k$, we have $\Ext^k(X_0,X_i)\cong
\Ext^{k-1}(X_0,X_{i-1})\cong \cdots \cong \Ext^{k-i}(X_0,X_0)=0$,
 since $0 < k-i < d+1$, and if $i > k$, we have  $\Ext^k(X_0,X_i)\cong
\Ext^{k+1}(X_0,X_{i+1})\cong \cdots \cong
\Ext^{k+d+1-i}(X_0,X_{d+1})=\Ext^{k+d+1-i}(X_0,X_0)=0$, since $0 <
k+d+1-i < d+1$.

\end{proof}

\begin{lem}

$\End X_i \cong \End X_0$ as algebras, for $0 \leq i \leq d$.

\end{lem}

\begin{proof}

We only need to prove the ring isomorphism End$X_1 \cong
\mbox{End}X_0$, since the others are done by induction. It is
exactly the same as the proof of the case $d=1$ in \cite{BMRRT}.

\end{proof}

\begin{lem}

$\dim_{\End X_i}\Ext^{k}(X_i,X_j) = \left\{ \begin{array}
{lll}1&\mbox{if } i+k-j=0 \text{ mod } (d+1) \\ 0 &  otherwise
\end{array} \right.$,  for $0 \leq k \leq d$. If we fix
 an $\End {X_i}-$basis $\{\delta_i\}$ of $\Ext^1(X_i,X_{i+1})$, then for
 any $0\le i\le d$ and $0\le k\le d$, $\Ext^{k}(X_i,X_{i+k})$ has
 an $\End (X_i)-$basis $\{\delta_{i+k}[k]\cdots \delta _{i+1}[1]\delta_i\}$,
 where $X_{i+k}=X_{i+k-(d+1)}$ and $\delta
 _{i+k}=\delta_{i+k-(d+1)}$,
 for $i+k>d$.

\end{lem}

\begin{proof}

The case of $i=0$ of the first part follows easily from the two
lemmas above, and the case for arbitrary $i$ follows from the same
proof after replacing $0$ by $i$. For the second part, it is easy to
see that any morphisms $\delta_{i+k}[k]\cdots \delta
_{i+1}[1]\delta_i$ are non-zero in  Ext$^{k}(X_i,X_{i+k})$, hence
 form a basis over End$X_i$ of  Ext$^{k}(X_i,X_{i+k})$.

\end{proof}

\begin{defn} A set of $d+1$ indecomposable objects $X_0,X_1,\cdots, X_d$ in
$\mathcal {C}_d(\mathcal{H})$ is called a exchange team if they
satisfy Lemma 5.3. i.e.  $\dim_{\End X_i}\Ext^{k}(X_i,X_j) \\
  =
\left\{
\begin{array} {lll}1&\mbox{if } i+k-j=0 \text{ mod } (d+1) \\ 0 & otherwise
\end{array} \right.,$
   for $0 \leq k \leq d$. If we fix an $\mbox{End}{X_i}-$basis $\{\delta_i\}$ of
 $\Ext^1(X_i,X_{i+1})$, then for any $0\le i\le d$ and $0\le k\le d$,
 $\Ext^{k}(X_i,X_{i+k})$ has an End$X_i-$basis $\{\delta_{i+k}[k]\cdots
\delta _{i+1}[1]\delta_i\}$, where $X_{i+k}=X_{i+k-(d+1)}$ and
$\delta _{i+k}=\delta_{i+k-(d+1)}$, for $i+k>d$.\end{defn}

This is a generalization of the notation of exchange pairs in
cluster categories, defined in \cite{BMRRT}.

\medskip

Given an exchange team $\{X_i\}^d_{i=0}$, by definition we can find
$d+1$ non-split triangles
\begin{displaymath} (***)\ \ X_{i+1}\s{g_i}{\longrightarrow}B_i\s{f_i}{\longrightarrow}X_i{\longrightarrow}X_{i+1}[1]\end{displaymath}
in $\mathcal {C}_d(\mathcal{H})$, where we use the same notation as
before. We will now start to prove that  $B=\bigoplus_{0 \leq i \leq
d}B_i$ is a rigid object.

\begin{lem}

With the notation above, we have
\begin{displaymath} \Ext^k(B \bigoplus X_i,B \bigoplus X_i)=0,\end{displaymath} for all $1 \leq k \leq
d$ and $0 \leq i \leq d$.

\end{lem}

\begin{proof}

Apply $\Hom(X_0,-)$ to the triangle
$X_1{\longrightarrow}B_0{\longrightarrow}X_0\s{\delta_0}{\longrightarrow}X_1[1]$
to get the exact sequence
\begin{displaymath}\Hom(X_0,X_0)\s{\alpha}{\longrightarrow}\Ext(X_0,X_1){\longrightarrow}
\Ext(X_0,B_0){\longrightarrow}\Ext(X_0,X_0).\end{displaymath}

Since $\alpha \neq 0$ ($\alpha(1_{X_0})=\delta_0 \neq 0$) and
dim$_{\End(X_0)}\Ext(X_0,X_1)=1$, while $\Ext(X_0,X_0)=0$ by
assumption, it follows that $\Ext(X_0,B_0)=0$. By assumption,
 $\Ext^{k}(X_0,X_1)=0$ and $\Ext^{k}(X_0,X_0)=0$ for any $2 \leq k
\leq d$, so it follows that $\Ext^{k}(X_0,B_0)=0$ for any $2 \leq k
\leq d$. Hence $\Ext^k(X_0,B_0)=0$, for $1 \leq k \leq d$.

Apply $\Hom(X_0,-)$ to the triangle
 $X_{i+1}\s{g_i}{\longrightarrow}B_i\s{f_i}{\longrightarrow}X_i{\longrightarrow}X_{i+1}[1]$
to get the exact sequence
\[ \begin{array}{*{3}{c@{\:{\s{\:}{\longrightarrow}}\:}}c}
& \Ext(X_0,X_{i+1}) & \Ext(X_0,B_i) & \Ext(X_0,X_i)\\
\multicolumn{4}{c}{\dotfill}\\
& \Ext^i(X_0,X_{i+1}) & \Ext^i(X_0,B_i) & \Ext^i(X_0,X_i)\\
& \Ext^{i+1}(X_0,X_{i+1}) & \Ext^{i+1}(X_0,B_i) & \Ext^{i+1}(X_0,X_i)\\
\multicolumn{4}{c}{\dotfill}\\
& \Ext^d(X_0,X_{i+1}) & \Ext^d(X_0,B_i) & \Ext^d(X_0,X_i).
\end{array} \]
 $\Ext^i(X_0,X_i){\longrightarrow}\Ext^{i+1}(X_0,X_{i+1})$ is
 an isomorphism (because $f \in \Ext^{i+1}(X_0,X_{i+1})$ can be
decomposed), and $\Ext^k(X_0,X_{i+1})=0=\Ext^l(X_0,X_i)$ for $k \neq
i+1$ and $l \neq i$, so $\Ext^{k}(X_0,B_i)=0$ for any $1 \leq k \leq
d$. Analogously, we get $\Ext^{k}(X_i,B_j)=0$ for all $1 \leq k \leq
d$ and $0 \leq i,j \leq d$.

Apply $\Hom(B,-)$ to the triangles
 $X_{i+1}{\longrightarrow}B_i{\longrightarrow}X_i{\longrightarrow}X_{i+1}[1]$
to get the exact sequences
\begin{displaymath}\Ext^k(B,X_{i+1}){\longrightarrow}\Ext^k(B,B_i){\longrightarrow}\Ext^k(B,X_i).\end{displaymath}
 Then $\Ext^k(B,B_i)=0$ for all $0 \leq i \leq d$ and $1 \leq k \leq d$,
so $\Ext^k(B,B)=0$ for all $1 \leq k \leq d$.

\end{proof}

Note that this implies that the $X_i$ cannot be direct summands of
$B$ (if $X_i \in \mbox{add}B$ for some $i$, then $\Ext(X_i,
X_{i+1})$ is a direct summand of $\Ext(B \bigoplus X_{i+1}, B
\bigoplus X_{i+1})=0$, a contradiction) and $B$ is a rigid object in
${\C_d(\mathcal{H})}$. Hence $B$ can be extended to a $d-$tilting
object by Corollary 3.4. Let $T=B \bigoplus T'$ be a $d-$cluster
tilting object in ${\C_d(\mathcal{H})}$.

\begin{lem}
 Under the same assumptions and notation as before, if $N$ is an
 indecomposable summand of $T$ and there exists some $j$ such that
$N$ is not isomorphic to $X_i$ for all $i \neq j$, then
$\Ext^{k}(N,X_j)=0$ for any $1 \leq k \leq d$.

\end{lem}

\begin{proof}

Assume by contradiction that $\Ext^{k}(N,X_j) \neq 0$ for some $1
\leq k \leq d$, and there is some indecomposable summand $N$ of $T$
with $N \ncong X_i$ for all $i \neq j$. Applying $\Hom(N,-)$ to
 the $d+1$ triangles $(***)$, we get $\Ext^{1}(N,X_{j-k+1}) \cong
\Ext^{k}(N,X_j) \neq 0$. Without loss of generality, we may assume
that $j-k=0$. So we have $\Hom(N,X_1[1]) = \Ext^1(N,X_1) \neq 0$ and
an exact sequence
\begin{displaymath}\Hom(N,X_0){\longrightarrow}\Hom(N,X_1[1]){\longrightarrow}0,\end{displaymath}
which implies that there exists a non-zero morphism $t \in
\Hom(N,X_0) \neq 0$ such that $\delta_0t \neq 0$. Applying
$\Hom(N,-)$ to the $d+1$ triangles $(***)$, we get
 $\Ext^{d}(N,X_{d}) \cong \Ext^{d-1}(N,X_{d-1})\cong \cdots
\Ext^1(N,X_1) \neq 0$, and then
$\delta_d[d]\cdots\delta_1[1]\delta_0t \neq 0$. Denote by
\begin{displaymath}X_0[d]{\longrightarrow}A\s{r}{\longrightarrow}X_0{\longrightarrow}X_0[d+1]\end{displaymath}
the AR-triangle ending at $X_0$ in $\C_d(\mathcal{H})$. Consider the
commutative diagram
$$ \begin{array}{ccclclcl}
X_0[d]&{\longrightarrow} &A&\s{r}{\longrightarrow}&X_0&\longrightarrow&X_0[d+1]\\
\parallel&&\downarrow b_1&&\downarrow b_2&& \parallel\\
X_0[d]&\s{g_d[d]}{\longrightarrow}&B_d[d]&\s{f_d[d]}{\longrightarrow}&X_d[d]&\s{\delta_d[d]}{\longrightarrow}&X_0[d+1],
\end{array}$$
where the map $b_1$ exists since $\delta _d[d] \neq 0$ (thus
 $g_d[d]$ is not a section), and hence there exists a map $b_2$ such that the
diagram commutes. From Definition 5.4, we know that
$\Hom(X_0,X_d[d])$ has an $\mbox{End}X_0-$basis
 $\{\delta_d[d]\cdots\delta_1[1]\delta_0\}$. Since $b_2 \in
\Hom(X_0,X_d[d])$ is not zero, there exists an isomorphism $\phi \in
\End(X_0)$ such that $b_2=\delta_d[d]\cdots\delta_1[1]\delta_0\phi$.
Let $s=\phi^{-1}t \in \Hom(N, X_0)$, then $b_2s \neq 0$. Since $N
\ncong X_0$, there is some map $s'$:$N{\longrightarrow}A$, such that
$s=rs'$. Note that $b_2s=b_2rs'=f_d[d]b_1s'$ is a non-zero map, and
consequently $b_1s' \neq 0$. But this contradicts
 $\Hom(N,B_d[d])=0$. This completes the proof of the lemma.

\end{proof}

\begin{lem}

If $\add(\bigoplus_{{1 \leq i \leq d}, {i \neq j}}X_i){\bigcap} \add
T=\{0\}$ for some ${1 \leq j \leq d}$, then $X_j$ is a direct
summand of $T$. Writing  $T$ as $X_j^k \bigoplus \overline{T}$,
where the $X_j$ are not direct summands of $\overline{T}$, then $X_i
\bigoplus \overline{T}$ is also a $d-$cluster tilting object for any
${0 \leq i \leq d}$.

\end{lem}

\begin{proof} The first assertion follows directly from Lemma 5.6.
The second follows from Theorem 4.3 and Lemma 5.6.

\end{proof}

 In summary, we have the following main result:

\begin{thm}
\label{mainthm2}

The $d+1$ rigid indecomposable objects $\{ X_i \}_{0 \leq i \leq d}$
form the set of complements of an almost complete d-cluster tilting
object in $\C_d(\mathcal{H})$ if and only if they form an exchange
team.

\end{thm}

Since the chain of $d+1-$triangles of the complements of an almost
complete $d-$cluster tilting object form a cycle, their distribution
is uniform. In particular there are two cases: either
 every complement has a different degree, or that the degree of any complement is smaller than $d-1$ and only two
complements have the same degree. We can summarize the cases
 as follows.

\begin{prop}

Suppose $\deg X_0=0$ and $\deg X_1 \neq 0$. Then there exists some
k, with $0 \leq k \leq d$, such that
$\deg X_i=\left\{\begin{array} {r@{\quad\mbox{if}\quad}l} d-i & 1 \leq i \leq k \\
d+1-i & k+1 \leq i \leq d . \end{array} \right.$

\end{prop}

\begin{proof}
By Lemma 3.1, we know that $\deg X_i \geq d-i$ for $1 \leq i \leq
d$. Since $d+1-$triangle chains form a cycle, analyzing the degree
in
 the opposite direction from $X_0$, we get $\deg X_i \leq d-i+1$
for $1 \leq i \leq d$. If $\deg X_1=d$, then $\deg X_2=d-1$, since
 $\Hom(X_1,X_2[1]) \neq 0$ forces $\deg X_2 \geq d-1$. By
induction, $\deg X_i=d-i+1$ for $1 \leq i \leq d$. This situation is
equivalent to $k = 0$. If $\deg X_1=d-1$, then there exists some $k$
 such that $\deg X_k = \deg X_{k+1}$. By the way of the case $\deg
X_1=d$, we obtain the conclusion.

\end{proof}

\section{Middle terms of the $d+1$ triangles}

Throughout this section, we assume that $d \geq 2$.   We assume that
$X$ is a basic almost complete $d-$cluster tilting object, and that
the $d+1$ complements $X_0,\cdots , X_d$ of $X$ are connected by
 the $d+1$ triangles in $(*)$ in Section 4:

\begin{displaymath}(*)\ \  X_{i+1}\s{g_i}{\longrightarrow}B_i\s{f_i}{\longrightarrow}
X_i\s{\delta_i}{\longrightarrow}X_{i+1}[1],\end{displaymath} where
 for $i=0, 1, \cdots, d$, $B_i \in \add X$, the map $f_i$ is the minimal right
 $\mbox{add}X-$approximation of $X_i$ and $g_i$ is the minimal left
add$X$-approximation of $X_{i+1}$. \vspace{10pt}

 In \cite{BMRRT}, there was a conjecture
that the sets of indecomposables of $B_i$ appeared in the triangles
$(*)$ are disjoint in cluster categories. That has been solved in
\cite{BMR}. We will prove the same statement for d-cluster
categories. Prior to this, we need some preparatory work. For a
tilting module $T$ in $\H$, any two non-isomorphic summands $T_1$,
$T_2$ of $T$ have the following property: $\Hom(T_1,T_2)=0$ or
$\Hom(T_2,T_1)=0$ (see \cite{Ker}). The same property holds for
$d-$cluster tilting objects in d-cluster categories when $d \geq 3$.

\begin{lem}

Suppose $d \geq 3$. Let $T_1$, $T_2$ be two non-isomorphic summands
of a $d-$cluster tilting object $T$ in $\C_d(\mathcal{H})$. Then
$\Hom(T_1,T_2)=0$ or $\Hom(T_2,T_1)=0$.

\end{lem}

\begin{proof}

If not, then $\Hom(T_1,T_2) \neq 0$ and $\Hom(T_2,T_1) \neq 0$. Then
$\deg T_1 = \deg T_2$ by the fact that $d \geq 3$ and Lemma 4.7 in
\cite{Zh3}. Let $k$ denote this common value. Then $T_1, \ T_2$ are
 of the forms $T_1'[k],\ T_2'[k]$ respectively, where $T'_1$ and $T'_2$
are partial tilting modules in $\mathcal{H}$. Hence $\Hom(T_1,T_2)
\cong \Hom_{\D}(T_1',T_2')\not= 0$ and $\Hom(T_2,T_1) \cong
\Hom_{\D}(T_2',T_1')\not=0$ \cite{Ker}. That is a contradiction.

\end{proof}

As a consequence, we get the following simple result.

\begin{lem}
Let $d \geq 3$. Then $\Hom(X_i,X_{i+1})=0$.
\end{lem}

\begin{proof}
Apply $\Hom(X_i,-)$ to the triangle
$X_{i+1}{\longrightarrow}B_i{\longrightarrow}X_i{\longrightarrow}X_{i+1}[1]$
to get the exact sequence
\begin{displaymath}\Hom(X_i,X_i[-1]){\longrightarrow}\Hom(X_i,X_{i+1}){\longrightarrow}\Hom(X_i,B_i).\end{displaymath}
In this exact sequence, $\Hom(X_i, X_i[-1])=0$ since $d\ge 3$. Since
$B_i{\longrightarrow}X_i$ is the minimal right add$X-$approximation,
 $\Hom(Y, X_i) \neq 0$ for any indecomposable direct summand $Y$ of
$B_i$. It follows from Lemma 6.1 that $\Hom(X_i,B_i)=0$. Thus
 $\Hom(X_i,X_{i+1})=0$.
\end{proof}

Now we are able to prove the main conclusion in this section.

\begin{thm}

Let  $\verb|{| B_i \verb|}|_{0 \leq i \leq d}$ be as above. Then the
sets of indecomposable summands of $B_i$, for $i=0, \cdots , d$,
 are disjoint.

\end{thm}

\begin{proof} We divide the proof into two cases:

(1). The case when $d=2$. Suppose $\deg X_0=0$. Assume by
contradiction that two of $B_0, B_1, B_2$ have non-trivial
intersection. Without loss of generality, we suppose that there
 exists an indecomposable object $T_1 \in \mbox{add}B_0 \bigcap
\mbox{add}B_1$. Then $\Hom(X_1,T_1) \neq 0 \neq \Hom(T_1,X_1)$,
 which implies that $\deg X_1 \neq \deg T_1$ (see \cite{Ker}). We
claim that $\deg X_1=1$, $\deg X_2=0$, and $\deg T_1=0$. If $\deg
X_1=0$, then $\deg T_1=0$ by Lemma 4.9 in \cite{Zh3}, a
contradiction. If $\deg X_1=2$ and $\deg T_1=0$, then
 $\Hom(T_1,X_1)=0$ by Lemma 4.7 in \cite{Zh3}, a contradiction. If
$\deg X_1=2$ and $\deg T_1=1$, then $\Hom(X_1, T_1)=0$ by Lemma 4.7
in \cite{Zh3}, a contradiction. So $\deg X_1=1$, and then $\deg
T_1=0$ (otherwise, $\deg T_1=2$ which implies $\Hom(T_1, X_1)=0$, a
contradiction). From Proposition 5.9, we have $\deg X_2=0$. Hence
 the degree of any indecomposable summands of $B_2$ is zero. Then
$\Hom(X_2,B_2)=0=\Hom(B_2,X_0)$ (see the discussion in the proof of
Lemma 6.2). Apply $\Hom(X_2,-)$ to the triangle
$X_0{\longrightarrow}B_2{\longrightarrow}X_2{\longrightarrow}X_0[1]$
to get the exact sequence
\begin{displaymath}\Hom(X_2,X_2[-1]){\longrightarrow}\Hom(X_2,X_0){\longrightarrow}\Hom(X_2,B_2),\end{displaymath}
where $\Hom(X_2,B_2)=0$, so $\Hom(X_2,X_0)=0$ (for any map $r \in
\Hom(X_2,X_0)$, there exists $s \in \Hom(X_2,X_2[-1])$ $\cong
\Hom(X_2,\tau^{-1}X_2[1])$ $\cong \Hom_{\D}(X_2,\tau^{-1}X_2[1])$
$\cong \Hom_{\D}(\t X_2[-2],X_2[-1])$  and  $t \in
\Hom(X_2[-1],X_0)$ $\cong \Hom(X_2,X_0[1])$ $\cong
\Hom_{\D}(X_2,X_0[1])$ $\cong \Hom_{\D}(X_2[-1],X_0)$ (both of the
second isomorphisms come from Lemma 4.2), such that $r=ts \in
\Hom_{\D}(\t X_2[-2],X_0)=0$). Write the second triangle in $(*)$ as
\begin{displaymath}X_2\s{\binom{h}{f}}{\longrightarrow}B_1' \bigoplus T_1
\s{\left( \alpha, \beta
\right)}{\longrightarrow}X_1{\longrightarrow}X_2[1],\end{displaymath}
where $\beta \in \Hom(T_1,X_1) \cong \Hom_{\D}(T_1,X_1)$.  Let $g$
be a non-zero map in $\Hom(T_1,X_0)$ (such a map exists because
$T_1$
 is a direct summand of $B_0$). Then we get ${\left( 0, g \right)}
\binom{h}{f}=gf \in \Hom(X_2,X_0)=0$, so there exists a map $\varphi
\in \Hom(X_1,X_0)$ $\cong \Hom(X_1, \t^{-1} X_0[2])$ $\cong
\Hom_{\D}(X_1, \t^{-1} X_0[2])$ (the second isomorphism come from
 Lemma 4.2) such that $\varphi (\alpha, \beta)=(0,g)$. Then
$g=\varphi \beta \in \Hom_{\D}(T_1, \t^{-1}X_0[2])=0.$ This is a
contradiction.
\medskip

(2). The case when $d \geq 3$. Suppose $T_1$ is an indecomposable
summand of both $B_i$ and $B_j$, $i<j$. Define
$d(B_i,B_j)=\mbox{min}\{ j-i, i-j+d+1 \}$.

If $d(B_i,B_j)=1$, then without loss of generality we may suppose
that $i=0$ and $j=1$; then $\Hom(X_1,T_1) \neq 0$ and $\Hom(T_1,X_1)
\neq 0$. But $X_1$ and $T_1$ are two non-isomorphic indecomposable
summands of a $d-$cluster tilting object $X_1 \bigoplus X$, which is
impossible by Lemma 6.1.

If $d(B_i,B_j)=2$, then without loss of generality we may suppose
that $i=1$ and $j=3$; then $\deg X_2 = \deg X_3 = \deg T_1$. Let $k$
 denote this common value. Then $\deg X_4 = k-1$ when $k \geq 1$, and
 $\deg X_4 = d-1$ when $k=0$. Apply $\Hom(X_2,-)$ to the triangle
$X_4\s{g_3}{\longrightarrow}B_3\s{f_3}{\longrightarrow}X_3\s{\delta_3}{\longrightarrow}X_4[1]$
to get an exact sequence
\begin{displaymath}\Hom(X_2,X_4){\longrightarrow}\Hom(X_2,B_3){\longrightarrow}\Hom(X_2,X_3).\end{displaymath}
Then $\Hom(X_2,X_4){\longrightarrow}\Hom(X_2,B_3)$ is an epimorphism
since $\Hom(X_2,X_3)=0$. Since $T_1 \in \add B_1$, there exists a
non-zero morphism $s \in \Hom(X_2,T_1)$, so the morphism
$\binom{s}{0}:X_2{\longrightarrow}T_1\bigoplus B_3'$ is not zero,
where $B_3=B_3' \bigoplus T_1$. Hence there exists $r \in
\Hom(X_2,X_4)$ such that $s=g_3r$. Let
$g_3=\binom{h}{h'}:X_4{\longrightarrow}T_1\bigoplus B_3'$, where $h
\in \Hom(X_4, T_1)$, then $s=hr$. Since $\Hom(X_2,X_4)\cong
\Hom_{\D}(X_2,\t^{-1}X_4[d])$ and $\Hom(X_4,T_1) \cong
\Hom_{\D}(\t^{-1}X_4[d], \t^{-1}T_1[d])$, it follows that $hr \in
\Hom_{\D}(X_2, \t^{-1}T_1[d])=0$, a contradiction.

If $d(B_i,B_j) \geq 3$, then the degrees of the summands of $B_i$
and $B_j$ are distinct. Hence the sets of indecomposable summands of
$B_i$ are disjoint, for $i=0, \ldots , d$.

\end{proof}

\section{Cluster combinatorics of $d-$cluster categories}

Denote by $\mathcal{E}(\mathcal{H})$ the set of isomorphism classes
of indecomposable rigid modules in $\mathcal{H}$. The set
$\mathcal{E}(\mathcal{C}_d(\mathcal{H}))$ of isoclasses of
indecomposable rigid objects in $\mathcal{C}_d(\mathcal{H})$ is the
(disjoint) union of the subsets $\mathcal{E}(\mathcal{H})[i]$,
$i=0,1, \cdots, d-1,$ with $\{ P_j[d]| 1\leq j\leq n\}$ (see Section
4 in [Zh3]). A subset $\mathcal{M}$ of
$\mathcal{E}(\mathcal{C}_d(\mathcal{H}))$ is called rigid if for any
$X, Y\in \mathcal{M}$, $\mathrm{Ext}^i(X,Y)=0$ for all $i=1,\cdots,
d$. Denote by $\mathcal{E}_+(\mathcal{C}_d(\mathcal{H}))$ the subset
of $\mathcal{E}(\mathcal{C}_d(\mathcal{H}))$ consisting of all
indecomposable exceptional objects other than $P_1[d], \cdots
P_n[d].$

Now we recall the definition of simplicial complexes associated to
the $d-$cluster category $\mathcal{C}_d(\mathcal{H})$ and the root
system $\Phi$ from [Zh3].

\begin{defn} The cluster complex $\Delta ^d(\mathcal{H})$ of $\mathcal{C}_d(\mathcal{H})$
is a simplicial complex with
 $\mathcal{E}(\mathcal{C}_d(\mathcal{H}))$ as its set of vertices,
 and the rigid subsets of $\mathcal{C}_d(\mathcal{H})$ as its
simplices. The positive part $\Delta ^d_{+}(\mathcal{H})$ is the
subcomplex of $\Delta ^d(\mathcal{H})$ on the subset
$\mathcal{E}_+(\mathcal{C}_d(\mathcal{H}))$.

\end{defn}

From the definition,  the facets (maximal simplices) are exactly
the $d-$cluster tilting subsets (i.e. the sets of indecomposable
objects of $\mathcal{C}_d(\mathcal{H})$ (up to isomorphism) whose
direct sum
 is a $d-$cluster tilting object).
\medskip

As consequences of results in Sections 3. 4. 5., we have that

\begin{prop}\begin{enumerate}
\item A face of the cluster complex $\Delta ^d(\mathcal{H})$ is a
facet if and only if it contains exactly $n$ vertices. In
particular, all facets in $\Delta ^d(\mathcal{H})$ are of size $n$
\item Every codimension $1$ face of $\Delta ^d(\mathcal{H})$ is
contained in exactly $d+1$ facets.
\item Any codimension $1$ face in $\Delta ^d(\mathcal{H})$ has complements of each color.
\end{enumerate}

\end{prop}

Throughout the rest of this section, we assume that $\mathcal{H}$
is the category of finite dimensional representations of a valued
quiver $(\G, \Omega, \mathcal{M})$. For basic material about
valued quivers and their representations, we refer to \cite{DR}.

Let $\Phi$ be the root system of the Kac-Moody Lie algebra
corresponding to the graph $\G$. We assume that $P_1, \cdots, P_n$
are the non-isomorphic indecomposable projective representations in
$\mathcal{H}$, and $E_1, \cdots, E_n$ are the simple representations
with dimension vectors $\alpha _1, \cdots, \alpha _n$, where
$\alpha_1, \cdots , \alpha_n$ are  the simple roots in $ \Phi.$ We
use $\Phi_{\geq -1}$ to denote the set of almost positive roots,
i.e. the set of positive roots together with the $-\alpha_i$.

\medskip

Fix a positive integer $d$, for any $\alpha \in \Phi ^+$, following
\cite{FR2}, we call $\alpha ^1, \cdots, \alpha^d$ the $d$
``colored'' copies of $\alpha$.

\begin{defn}\cite{FR2} The set of colored almost positive roots is
$$\Phi^d_{\geq-1}=\{ \alpha ^i: \alpha\in \Phi_{>0}, i\in \{1, \cdots, d\}\}\bigcup
\{(-\alpha _i)^1: 1\leq i\leq n\ \}.$$

\end{defn}

\medskip

We now define a map $\gamma ^d_{\mathcal{H}}$ from $\mathrm{ind}
\mathcal{C}_d(\mathcal{H})$ to $\Phi ^d_{\ge-1} $. Note that any
indecomposable object $X$ of degree $i$ in
$\mathcal{C}_d(\mathcal{H})$ has the form $M[i]$, for some $M\in
\mathrm{ind}\mathcal{H}$, and if $i=d$ then $M=P_j$, an
indecomposable projective representation.

\begin{defn}

Let  $\gamma^d_{\mathcal{H}}$ be defined as follows. Let $M[i]\in
\mathrm{ind}\mathcal{C}_d(\mathcal{H})$, where $M\in
\mathrm{ind}H$ and $i\in \{ 1, \cdots, d\}$ (note that if $i=d$
then $M=P_j$ for some $j$). We set
$$\gamma^d_{\mathcal{H}}(M[i])=\left\{
\begin{array}{lrl}(\underline{\mathrm{dim}}M)^{i+1} & \mbox{ if } & 0\leq i\leq d-1 ;\\
&&\\
(-\alpha _j)^1& \mbox{ if } &i=d,\end{array}\right.$$

\end{defn}

Note that if $\G$ is a Dynkin diagram, then $\gamma^d_{\mathcal{H}}$
is a bijection.

We denote by $\Phi^{sr}_{>0}$ the set of real Schur roots of $(\G,
\Omega)$, i.e.
$$\Phi^{sr}_{>0}=\{ \underline{\mathrm{dim}} M: \ M\in \mathrm{ind} \mathcal{E(H)}\ \}.$$
Then the map $M\mapsto \underline{\mathrm{dim}} M$ gives a 1-1
correspondence between $\mathcal{E(H)}$ and $\Phi^{sr}_{>0}$
\cite{Rin}.

If we denote the set of colored almost positive real Schur roots by
$\Phi^{sr, d}_{\geq -1}$ (which consists by definition of $d$ copies
of the set $\Phi^{sr}_{>0}$ together with one copy of the negative
simple roots), then the map   $\gamma^d_{\mathcal{H}}$ gives a
bijection from $\mathcal{E}(\mathcal{C}_d(\mathcal{H}))$ to
$\Phi^{sr, d}_{\geq -1}$.  $\Phi^{sr, d}_{\geq -1}$ contains a
subset $\Phi^{sr,d}_{>0}$ consisting of all colored positive real
Schur roots. The restriction of $\gamma^d_{\mathcal{H}}$ gives a
bijection from $\mathcal{E}_+(\mathcal{C}_d(\mathcal{H}))$ to
$\Phi^{sr, d}_{>0}$.

Using this bijection, in \cite{Zh3} we defined, for any root system
$\Phi$ and $\mathcal{H}$, an associated simplicial complex $\Delta
^{d,\mathcal{H}}(\Phi)$ on the set $\Phi^{sr, d}_{>0}$, which is
called the generalized cluster complex of $\Phi$ and is a
generalization of the generalized cluster complexes defined by
Fomin-Reading \cite{FR2}, see also \cite{Th} for finite root systems
$\Phi$. It was proved that $\gamma^d_{\mathcal{H}}$ defines an
isomorphism from the simplicial complex $\Delta ^d(\mathcal{H})$ to
the generalized cluster complex $\Delta^{d,\mathcal{H}}(\Phi),$
which sends vertices to vertices, and $k-$faces to $k-$faces
\cite{Zh3}.

\medskip

\begin{cor}
\begin{enumerate}
\item A face of the generalized cluster complex
$\Delta ^{d,\mathcal{H}}(\Phi)$ is a facet if and only it contains
exactly $n$ vertices. In particular, $\Delta ^{d,\mathcal{H}}(\Phi)$
is  of pure dimension $n-1$.

\item Any codimension $1$ face of $\Delta ^{d,\mathcal{H}}(\Phi)$ is contained in exactly $d+1$ facets.

\item For any codimension $1$ face of $\Delta
^{d,\mathcal{H}}(\Phi)$, there are complements of each color.

\end{enumerate}
\end{cor}

\begin{proof} Combining Proposition 7.2. with the fact that  $\gamma^d_{\mathcal{H}}$ is an isomorphism from  $\Delta ^d(\mathcal{H})$ to
 $\Delta^{d,\mathcal{H}}(\Phi)$ [Zh3], we have all the conclusions in the corollary.

\end{proof}

\begin{center}
\textbf {ACKNOWLEDGMENTS.}\end{center} The authors would like to
thank Idun Reiten for her interest in this work. After completing
this work, the second author was informed by Idun Reiten that Anette
Wraalsen also proved Theorem 4.3 in \cite{Wr}; he is grateful to
Idun Reiten for this!

The authors would like to thank the referee for his/her very useful
suggestions to improve the paper.

\end{document}